\documentclass[12pt,psamsfonts]{article}

\usepackage{amsmath,amsfonts,amssymb,amsmath,amscd}
\usepackage[active]{srcltx}
\usepackage{graphicx,array}
\usepackage{amsthm}
\usepackage{amsbsy}
\usepackage{tikz}
\usepackage{float}

\newcolumntype{L}{>{\displaystyle}l}
\newcolumntype{C}{>{\displaystyle}c}
\newcolumntype{R}{>{\displaystyle}r}

\newcommand{\dem}{\noindent{\it Proof:\ }}
\newcommand{\eproof}{\noindent\mbox{\framebox [0.6ex]{}}  \medskip}
\newcommand{\demm}{\noindent{\it Proof of the Theorem \ref{Teorema:principal}:\ }}

\def \R{\mbox{${\mathbb R}$}} 
\def \C{\mbox{${\mathbb C}$}} 
\def \N{\mbox{${\mathbb N}$}} 
\def \s{\mbox{\rm{${\textbf{S}}$}}}
\def \r{\mbox{\rm{${\textbf{R}}$}}}

\def \z{\mbox{\rm{${\textbf{Z}}$}}}

\def \P{\mbox{${\mathcal{P}}$}}
\def \PP{\mbox{${\vec{\mathcal{P}}}$}}
\def \Q{\mbox{${\mathcal{Q}}$}}        
\def \QQ{\mbox{${\vec{\mathcal{Q}}}$}} 
\def \RR{\mbox{${\vec{R}}$}}
\def \SS{\mbox{${\vec{S}}$}}

\newcommand{\adl}{Ad_L}
\newcommand{\adltk}{Ad_{L^t}^{\:k}}
\newcommand{\adlk}{Ad_{L}^{\:k}}

\DeclareMathOperator{\im}{Im}
\DeclareMathOperator{\aut}{Aut}
\DeclareMathOperator{\re}{Re}

\newtheorem{theorem}{Theorem}[section]
\newtheorem{lemma}[theorem]{Lemma}
\newtheorem{proposition}[theorem]{Proposition}

\newtheorem{remark}[theorem]{Remark}

\begin{document}

\title{Normal form theory for  reversible equivariant vector fields}


\author{P. H. Baptistelli\\
{\small Department of Mathematics, UEM}\\
{\small Av. Colombo, 5790, 87020-900, Maring\'a -  PR, Brazil (phbaptistelli@uem.br)}
\and M. Manoel\\
{\small Department of Mathematics, ICMC - USP}\\
{\small C.P. 668, 13560-970 S\~ao Carlos - SP, Brazil (miriam@icmc.usp.br)}
\and I.O. Zeli\\
{\small Department of Mathematics, IMECC - UNICAMP}\\
{\small C.P. 13081-970, Campinas -  S\~ao Paulo, Brazil (irisfalkoliv@ime.unicamp.br)}
}
\date{}
\maketitle

\noindent {\it MSC}:  7C80, 34C20, 13A50\\

\noindent {\it Keywords}: {Normal form, reversibility, symmetry, homological operator \\


\begin{abstract}
We give a method to obtain formal normal forms of reversible equivariant vector fields. The procedure we present is based on the classical
 method of normal forms combined with tools from invariant theory. Normal forms of two classes of resonant cases are presented, both with linearization  having a $2$-dimensional nilpotent part and a semisimple part with purely imaginary eigenvalues.
 \end{abstract}

\section{Introduction}

Normal form theory has been used as a tool for the local study of the qualitative behavior of vector fields. The subject has been developed for many years, since Poincar\'e \cite{poincare}, Birkhoff \cite{birkhoff} and Dulac \cite{dulac}, Belitskii \cite{belitskii} and Takens \cite{takens}.  It is useful to have normal form procedure for vector fields with additional structures. Examples are equivariant vector fields, symplectic, volume preserving, reversible and combinations thereof. Many authors have used normal form theory in distinct contexts to study  limit cycles, family of periodic orbits, relative equilibria and relative periodic solutions  (see, for example \cite{marco1, melbourne, lima, marco2, mereu}).
The resulting normal forms in each of these contexts can then be expressed in terms of the  group
invariants. That, however, is independent of the normal form procedure.

The classical method consists of  performing changes of coordinates  around a singular point that are perturbations of the identity,
 $\xi = I + \xi_k,$ where $k \geq 2$ and $\xi_k$ is a homogeneous polynomial of degree $k.$ The aim
 is to annihilate  as many terms of degree $k$ as possible in the original vector field, obtaining a conjugate vector field
 written in a simpler and more convenient form. The method developed by Belitskii \cite{belitskii}  reduces
 this problem to computing the kernel of the so-called homological operator. The operator corresponds
 to the Lie bracket of the Lie algebra $\mathcal{X}$ of smooth vector fields  with a critical point at the origin. Now $\mathcal{X}$ can be made into a filtered Lie algebra $\mathcal{X}$ =
 $\mathcal{X}_1 \supset \mathcal{X}_2 \supset \cdots, $ where $\mathcal{X}_k$ is the set of
 vector fields of degree $k$ or higher. Let $X \in \mathcal{X}$ be a vector field, which we write as
 $X = X_1 + X_2 + \cdots$, where $X_k \in \mathcal{X}_k$ is homogeneous of degree $k$.
 The method consists of transforming $X$ at each level $k$ into normal form by successive transformations. On each level the method consists of restricting the homological operator to the  vector space of homogeneous polynomials of degree $k$, which is  associated to the adjoint $X_1^t$ of the linearization $X_1$ of the original vector field $X$. The normal form procedure preserves the filtration. Also, the normal form automatically belongs to the same Lie algebra as the original
 vector field when the transformations are restricted accordingly. In the present paper, the
 vector field is reversible equivariant under the action of a group $\Gamma$ and
the normal form inherits  the symmetries and reversing symmetries if the changes of coordinates are equivariant under $\Gamma.$

Elphick {\it et al.} \cite{elphick} give an algebraic method to obtain the
 normal form developed by Belitskii by choosing nonlinear terms that are equivariant under a one-parameter group $\s$ given by
\begin{equation}
\s=\overline{\left\{ e^{ sX_1^t}, s\in \R \right\}}.
\label{eq:grupo-s}
\end{equation}
When $X$ is $\Gamma-$equivariant (all elements of $\Gamma$ are symmetries),  the truncated normal
  form is $\s \times \Gamma-$equivariant (see \cite[Theorem XVI 5.9]{golub}). In this paper, we prove that  if $X$ is $\Gamma-$reversible-equivariant, then the truncated normal
  form is $\s \rtimes \Gamma-$reversible-equivariant (Theorem~\ref{Teorema:principal2}). We then use algebraic tools from invariant theory
  to apply results of Antoneli {\it et al.} \cite{manoel1} that will produce the normal form truncated at any degree. Note
  that in general the group $\s$  is noncompact, and compacteness  is assumed in the algorithms in \cite{manoel1}.  However, we observe that these
  can equally be used whenever the ring of invariants and the module of equivariants by $\s$ are finitely generated. \\

We have organized this paper as follows. In Section \ref{section:PRELIMINARIES} we  briefly introduce the notation and basic concepts about reversible equivariant theory. In Section \ref{section:SEMIDIRECT} we present results in invariant theory of the semidirect product of two arbitrary groups from the invariant theory of each group separately. These results are the grounds for normal form theory that is developed in Section~\ref{section:BNF}. Section~\ref{section:BNF} is devoted to a short introduction to the basic concepts of normal form theory and, from that, we introduce the reversible equivariant normal form and prove the main results, Theorem~\ref{Teorema:principal} and Theorem~\ref{Teorema:principal2}. Finally, in Section~\ref{section:APPLICATIONS} we apply the results of Section~\ref{section:BNF} to obtain the normal form for two types of vector fields  under $\z_2$-action and $\z_2 \times \z_2$-action. We consider a vector field whose  linearization has a $2$-dimensional
nilpotent part and two resonant purely imaginary eigenvalues. The first case, namely under the action of $\z_2,$ generalizes some normal forms with resonance presented by Lima and Teixeira in \cite{lima}.  A complete classification of nonresonant and resonant vector fields under the action of a group generated by two involutions appears in \cite{iris}.


\section{Preliminaries on invariant theory}

\label{section:PRELIMINARIES}

Let $\Gamma$ be a compact Lie group  acting linearly on a finite-dimensional real vector space $V$. Consider an epimorphism
\begin{equation}\label{defisigma}
\sigma: \Gamma \to \z_2 =\{\pm 1 \},
\end{equation}
where $\sigma(\gamma)= 1$ if $\gamma$ is a symmetry and $\sigma(\gamma)=-1$ if $\gamma$ is a reversing symmetry. We denote by $\Gamma_+$ the group of symmetries of $\Gamma$. So $\Gamma_+ = \ker \sigma$ is a normal subgroup of $\Gamma$ of index $2,$ and $\Gamma = \Gamma_+ \dot{\cup} ~\delta \Gamma_+$ for an arbitrary reversing symmetry $\delta \in \Gamma$. Throughout this work, $\delta$ denotes a reversing symmetry whenever it appears in the text.

A polynomial function $f: V \to \R$ is $\Gamma-$invariant if
\[f(\gamma x)=f(x), ~\forall \gamma \in \Gamma, \ x \in V,
\]and it is called $\Gamma-$anti-invariant if
\[
f(\gamma x)=\sigma(\gamma)f(x), ~\forall \gamma \in \Gamma, \ x \in V.
\]

We denote by $\P_V(\Gamma)$ the ring of the $\Gamma-$invariant polynomial functions and by $\Q_V(\Gamma)$ the module of the $\Gamma-$anti-invariant polynomial functions over the ring $\P_V(\Gamma)$.

A polynomial mapping $g: V \to V$ is $\Gamma-$equivariant if
\begin{equation}\nonumber
g(\gamma x)=\gamma g(x),~\forall \gamma \in \Gamma, \ x \in V,
\end{equation}
and it is  $\Gamma-$reversible-equivariant
\begin{equation}\nonumber
g(\gamma x)=\sigma(\gamma)\gamma g(x),~\forall \gamma \in \Gamma, \ x \in V.
\end{equation}

We denote by $\PP_V(\Gamma)$ the module of the $\Gamma-$equivariant polynomial mappings and by $\QQ_{V}(\Gamma)$ the module of the $\Gamma-$reversible-equivariant polynomial mappings, both over the ring $\P_V(\Gamma)$. The modules $\Q_V(\Gamma),$ $\PP_{V}(\Gamma)$ and $\QQ_V(\Gamma)$ are finitely generated and graded over the ring $\P_V(\Gamma),$ which is also finitely generated and graded (see \cite{manoel1}). When $\sigma$ is trivial, then $\P_V(\Gamma)$ and $\Q_V(\Gamma)$, as well as $\PP_V(\Gamma)$ and $\QQ_V(\Gamma)$, coincide.

Consider now the space $\P_V$ of polynomial functions $V \to \R$ and the space $\PP_V$ of polynomial mappings $V \to V$. Consider also the actions of $\Gamma$ on these spaces induced by the action of $\Gamma$ on $V$:
\begin{equation}
\label{eq:action-on-pol.-spaces}
\begin{tabular}{ccccc}
$\Gamma \times \P_V \to \P_V $  &  &   \text{and}    &   & $\Gamma \times \PP_V \to \PP_V ,$ \\
$~(\gamma, f) \mapsto \gamma \odot f $  & &&      & $~(\gamma, g) \mapsto \gamma \star g$
\end{tabular}
\end{equation}
where $\gamma \odot f (x) = f(\gamma x)$ and $\gamma \star g (x) = \gamma^{-1} g(\gamma x)$, $\forall x \in V$, $\forall \gamma \in \Gamma$. We then define the Reynolds operators on $\P_V(\Gamma_+)$, namely  $R, S : \P(\Gamma_+)\to \P(\Gamma_+)$, by
$$
R(f)=\frac{1}{2} \sum_{\gamma \Gamma_+ \in \Gamma/\Gamma_+} \gamma  \odot  f = \frac{1}{2}\left(f + \delta \odot f \right),
$$
$$
S(f) =\frac{1}{2} \sum_{\gamma \Gamma_{+} \in \Gamma/\Gamma_+} \sigma(\gamma) \gamma \odot f = \frac{1}{2}\left(f - \delta \odot f \right),
$$
and their equivariant versions  $\RR, \SS:\PP_{V}(\Gamma_+) \to
\PP_{V}(\Gamma_+),$
$$
\RR(g)=\frac{1}{2}  \sum_{\gamma \Gamma_{+} \in \Gamma/\Gamma_+} \gamma \star g = \frac{1}{2}\left(g + \delta \star g \right),
$$

$$
\SS(g)= \frac{1}{2}  \sum_{\;\gamma \Gamma_{+} \in \Gamma/\Gamma_+} \sigma(\gamma)\gamma \star g = \frac{1}{2}\left(g - \delta \star g \right).
$$

In \cite{manoel1}, these  have been used to prove the decompositions of $\P_V(\Gamma)-$modules
$$
\P_V(\Gamma_+)=\P_V(\Gamma) \oplus \Q_V(\Gamma) \quad \text{and} \quad \PP_V(\Gamma_+)=\PP_V(\Gamma) \oplus \QQ_V(\Gamma),
$$

\noindent which are used to give algorithms to compute generators of $\Q_V(\Gamma)$ and $\QQ_V(\Gamma)$ from the knowledge of generators of $\P_V(\Gamma_+)$ and $\PP_V(\Gamma_+)$. These algorithms are applied here, in Section \ref{section:APPLICATIONS}.

A result in \cite{manoel2} provides a simple way to compute a set of generators of  $\P_V(\Gamma)$ from generators of $\P_V(\Gamma_+)$. The result is stated below:

\begin{theorem}\cite[Theorem 3.1]{manoel2} \label{teoremagamma} Let $\Gamma$ be a compact Lie group acting on $V$ and let $\sigma:\Gamma \to \z_2$ an epimorphism as in $(\ref{defisigma}).$ Let $u_1, \ldots, u_s$  be a Hilbert basis for the ring $\P_V(\Gamma_+)$. Then the set
\begin{equation}\nonumber
\left\{ R(u_i), S(u_i)S(u_j), \:\forall \ 1 \leq i,j \leq s  \right\}
\end{equation}
is a Hilbert basis for the ring $\P_V(\Gamma).$
\end{theorem}


\section{Invariant theory for cartesian product of groups} \label{section:SEMIDIRECT}

Given two groups $\Gamma_1$ and $\Gamma_2$, recall that a semidirect product $\Gamma_1 \rtimes \Gamma_2$ is the direct product  $\Gamma_1 \times \Gamma_2$ as a set with a group operation induced by a homomorphism $\mu : \Gamma_2 \to \aut (\Gamma_1)$. If $\mu$ is trivial, the groups commute
and semidirect is in fact direct product as a group. We consider $\Gamma_1$ and $\Gamma_2$ acting on $V$ and denote $(\rho, V)$ and $(\eta, V)$
 their representations, respectively. From these, define the operation $(\Gamma_1 \rtimes \Gamma_2) \times V \to V,$ \begin{equation}
(\gamma_1,\gamma_2) v = \rho(\gamma_ 1) (\eta(\gamma_2) v).
\label{acaoproduto}
\end{equation}

We then have:

\begin{proposition} \label{prop:representacao} The operation $(\ref{acaoproduto})$ defines an action of the semidirect product $\Gamma_1\rtimes \Gamma_2$ on $V$ if, and only if,  the representation of $\mu(\gamma_2)(\gamma_1)$ is a conjugation, that is, \[\rho(\mu(\gamma_2)(\gamma_1))=\eta(\gamma_2)\rho(\gamma_1)\eta(\gamma_2)^{-1}.\]
\end{proposition}

The proof of the proposition above is direct from the definition of an action.

We also have that
\begin{equation}
\label{eq:acaoproduto}
\rho(\gamma_ 1)\eta(\gamma_2)  = \eta(\gamma_ 2) \rho(\mu(\gamma_2^{-1})(\gamma_1)),
\end{equation}
which highlights the non-commutativity of the $\Gamma_1$ and $\Gamma_2$ actions if and only if $\mu$ is nontrivial.
Yet, we notice that when $\Gamma_1$ has both symmetries and reversing symmetries, then for each $\gamma_2 \in \Gamma_2,$ the automorphism $\mu(\gamma_2)$ preserves the symmetries and reversing symmetries of $\Gamma_1$. \\

In this work, we assume that $\Gamma_1$ and $\Gamma_2$ admit a semidirect product with a representation in the conditions of Proposition \ref{prop:representacao}. We remark that this assumption may fail: consider the groups $\Gamma_1= \left\langle  \kappa_1  \right\rangle$ and $\Gamma_2= \left\langle  \kappa_2  \right\rangle$ acting on $\R^2$ as $\kappa_1(x,y)=(y,x)$ and $\kappa_2(x,y)=(x,-y)$. In this case, there is no representation of the direct product $\Gamma_1 \times \Gamma_2$ on $\R^2$ as given in the Proposition \ref{prop:representacao}. Nevertheless, there may be applications where we wish to work with $\Gamma = \left\langle \kappa_1, \kappa_2  \right\rangle$. Obviously, in this case, we have $\Gamma = \tilde \Gamma_1 \rtimes \tilde \Gamma_2$, with $\tilde \Gamma_1= \langle \kappa_1 \kappa_2  \rangle$ and $\tilde\Gamma_2 = \langle \kappa_2 \rangle$, and apply our theory to these two new groups.

If we now consider  $\Gamma_2$ endowed with an epimorphism $\sigma$ as in (\ref{defisigma}) and $\Gamma_1$ as a group of symmetries, we can naturally
 introduce an epimorphism on $\Gamma_1 \rtimes \Gamma_2 $ in order to preserve this structure:
\begin{equation}\label{eqsigmatil}
\begin{tabular}{cccl}
$\tilde{\sigma}:$ &$ \Gamma_1 \rtimes \Gamma_2 $ & $\rightarrow$ & $\z_2$ \\
          &$(\gamma_1, \gamma_2)$        &$\mapsto $     &$ \sigma(\gamma_2)$.
\end{tabular}
\end{equation}
The next result organizes the invariant theory of $\Gamma_1 \rtimes \Gamma_ 2$ in this setting:

\begin{proposition} \label{teoremainter} Let $\Gamma_1$ and  $\Gamma_2$ be compact Lie groups acting linearly on $V$ and consider the epimorphism $\tilde{\sigma}$ defined by $(\ref{eqsigmatil}).$ Then:
\begin{enumerate}
    \item [(i)\:] $\P_V(\Gamma_1 \rtimes \Gamma_2)= \P_V(\Gamma_1) \cap \P_V(\Gamma_2)$;
    \item [(ii)\:] $\PP_V(\Gamma_1 \rtimes \Gamma_2)= \PP_V(\Gamma_1) \cap \PP_V(\Gamma_2)$;
    \item [(iii)\:] $\Q_V(\Gamma_1 \rtimes \Gamma_2)= \P_V(\Gamma_1) \cap \Q_V(\Gamma_2)$;
    \item [(iv)\:] $\QQ_V(\Gamma_1 \rtimes \Gamma_2)= \PP_V(\Gamma_1) \cap \QQ_V(\Gamma_2)$.
\end{enumerate}
\end{proposition}
\dem
For $f \in \P_V(\Gamma_1 \rtimes \Gamma_2)$, take $\gamma_1=I$ to get that $f(\eta(\gamma_2) v)= f(v)$ for all $\gamma_2 \in \Gamma_2$, so $f\in
\P_V(\Gamma_2).$ Analogously, for $\gamma_2=I$ it follows that $f\in
\P_V(\Gamma_1).$ On the other hand, if $f \in \P_V(\Gamma_1) \cap \P_V(\Gamma_2)$, then $f((\gamma_1,\gamma_2)  v)=f(\rho(\gamma_1) (\eta(\gamma_2)  v)) = f(\eta(\gamma_2)  v) = f(v)$. The proofs of the other three equalities are similar to this.
\eproof

The proposition above is applied in Section \ref{section:BNF}
where the algebraic approach is applied to $\s \rtimes \Gamma$ for an appropriate group $\s$ of symmetries. \\

We end this section recalling the diagonal representation of $\Lambda_1 \times \Lambda_2$ on $V \times W$  of two commuting groups  $\Lambda_1$
and $\Lambda_2$ acting respectively on vector spaces $V$ and $W$:
$$
(\lambda_1, \lambda_2)(v,w) = ( \lambda_1 v, \lambda_2 w).
$$
It is immediate that the relations given in Proposition \ref{teoremainter} hold equally for the direct product. We also have:

\begin{lemma} \label{lemma:geradores} Let $\Lambda_1$ and $\Lambda_2$ be commuting groups acting linearly on $V$ and $W$, respectively. Let $\{u_1, \ldots, u_r\}$ and $\{\alpha_1, \ldots, \alpha_s\}$ be a Hilbert bases for $\P_V(\Lambda_1)$ and $\P_W(\Lambda_2),$ respectively. Let $\left\{ f_1, \ldots f_m \right\}$ be a generating set of $\PP_V(\Lambda_1)$ over the ring $\P_V(\Lambda_1)$ and let $\left\{ g_1, \ldots g_n \right\}$ be a generating set of $\PP_W(\Lambda_2)$ over the ring $\P_W(\Lambda_2)$. Then, $\{u_1, \ldots u_r, \alpha_1, \ldots, \alpha_s\}$ is a Hilbert basis for $\P_{V\times W}(\Lambda_1 \times \Lambda_2)$ and
$$
\left\{ \begin{pmatrix} f_i \\ 0_W \end{pmatrix}, \begin{pmatrix} 0_V \\ g_j \end{pmatrix}, \ 1 \leq i \leq m, \ \:1 \leq j \leq n\right\}
$$
is a generating set of the module $\PP_{V\times W}(\Lambda_1 \times \Lambda_2)$ over the ring $\P_{V\times W}(\Lambda_1 \times \Lambda_2)$, where $0_V$ and $0_W$ denote the zero vectors of $V$ and $W$, respectively.
\end{lemma}
\dem
The results follow directly from
\[\P_{V\times W}(\Lambda_1 \times \Lambda_2) = \P_{V\times W}(\Lambda_1 \times {\bf 1}) \cap
\P_{V\times W}({\bf 1} \times \Lambda_2) \]
and from the corresponding equality for the module of equivariants.

\section{Computing the normal form by  algebraic methods} \label{section:BNF}

In this section we provide an algebraic method to obtain the normal form of $\Gamma-$reversible equivariant vector fields on a finite dimensional vector space $V$. As we shall see, the method takes the action of $\Gamma$ into account to simplify the deduction of the normal form significantly. We now follow \cite{chow} to briefly recall the normal form theory without symmetries. Our approach adapts  that method.

Consider a system of ODEs
\begin{equation}
\label{eq:sistemainicial}
\dot{x}= X(x),~ x\in V.
\end{equation}
The interest of the dynamics analysis in normal form theory is local, around a singular point that we
assume to be the
origin, so $X(0) = 0$;  $X$ is assumed to be smooth around that point.  However, we focus now on a purely algebraic derivation of
the normal form (\ref{eq:sistemainicial}), which constitutes the first step in the dynamical analysis. Hence, we simplify the notation and use $V$ for the domain throughout.

From now on $L$ denotes the linear part of $X$ at the origin and consider the Taylor expansion of $X$ about that point
\begin{equation}
\label{eq:forma-de-h}
L(x) + X_2(x) + X_3(x) + \ldots \,
\end{equation}
with $X_k$ homogeneous of degree $k$, $k \geq 2$.  The method to obtain the normal form for the system (\ref{eq:sistemainicial}) consists of successive changes of coordinates on the source of the form $I + \xi_k,$ for $k\geq2$, where $I$ is the identity and $\xi_k$ is a homogeneous polynomial of degree $k$.  As a result, the new system has a ``simpler" form at its degree-$k$ level, no effect on its lower-order terms and with any form of its terms of degree greater than $k$. Here we do not deal necessarily with analytic mappings, so the systems are formally conjugate in the sense that their corresponding formal Taylor series are conjugate as formal vector fields.

After the changes of coordinates up to degree $k$ and rewriting in the $x$ variable, one obtains the new system with an intermediate form
\begin{equation} \label{sistemanovo}
\dot{x}=Lx + \sum_{n=2}^{k-1} g_n(x) +\tilde{X}_k(x)-((D\xi_k)_{(x)}Lx -L\xi_k(x)) + O|x|^{k+1}.
\end{equation}

\noindent Now consider the homological operator $\adl: \PP_V \to \PP_V$ given by
\begin{equation}
\adl(p)(x) = (Dp)_{(x)}Lx-Lp(x),
\end{equation}
where $\PP_V$ is the vector space of polynomial mappings $V \to V$. Since $\PP_V$ is a graded algebra, that is $\PP_V= \bigoplus^{\infty}_{k=0} \PP_V^{k}$, we can consider for each $k \geq 0$ the operator $\adlk,$ the restriction of $\adl$ to $\PP_V^{k}$. From (\ref{sistemanovo}) it follows that the new system can be simpler if its term of degree $k$ is of the form $g_k=\tilde{X}_k - \adlk(\xi_k)$, for some  $\xi_k$. This is clearly not unique, and a choice for $g_k$ was proposed by Belitskii \cite{belitskii} through an appropriate complement, $(\im\adlk)^c = \ker \adltk$, where $L^t$ denotes the adjoint operator of $L$. So, the method consists of determining the polynomials in the kernel of $\adltk$, and (\ref{eq:sistemainicial}) turns out to be formally conjugate to
\begin{equation} \nonumber
\dot{x}= Lx + g_2(x) + g_3(x) + \ldots
\end{equation}
with $g_k\in \ker \adltk,~k \geq 2$.

An alternative method to find the complement has been proposed by Elphick {\it et al.} in \cite{elphick} recognizing a group of symmetries on this space. More precisely, from the linear part $L$, consider the group $\s$ defined as in (\ref{eq:grupo-s}). So, $\s$ is a one-parameter closed subgroup of $GL(n)$ acting on $V$ by matrix product. The authors prove that
$\ker \adltk =\PP_V^k(\s)$ (see \cite[Theorem 2]{elphick}) and, therefore,
\begin{equation}
\label{eq:ELPHICK}
\PP_V^k= \PP_V^k(\s) \oplus \adl(\PP_V^k).
\end{equation}

Now, we are interested in obtaining the normal form for a system of the form (\ref{eq:sistemainicial}) when it is $\Gamma-$reversible-equivariant. Changes of coordinates in this setting are assumed to be $\Gamma-$equivariant, so that all the symmetries and reversing symmetries of the original system are preserved in the normal form. This could be done by first finding $\ker \adltk$ and after that imposing the existence of symmetries and reversing symmetries in the normal form. The result that we show here provides a different method for this process, based on the knowledge of the invariant theory for $\Gamma$ and $\s$, as an extension of the method presented in \cite[Chapter XVI, $\S 5$]{golub}. Although there are applications for which the group $\s$ is not compact, the theory developed in \cite{manoel1} can be applied here, as long as there is a finite set of generators for the ring $\P_V(\s)$ and for the module $\PP_V(\s)$ over $\P_V(\s)$.

The next two lemmas are useful in the remainder of this section.

\begin{lemma}  The homological operator $\adl$ interchanges the modules in the sum decomposition $\PP_V(\Gamma_+)=\PP_V(\Gamma) \oplus \QQ_V(\Gamma)$.
\end{lemma}
\dem If $p \in \PP_V(\Gamma),$ then
\begin{align}
  \adl (p)(\gamma x) &= (Dp)_{(\gamma x)}L(\gamma x) - L p(\gamma x) \nonumber\\
                      &= \gamma (Dp)_{(x)}\gamma^{-1} \sigma(\gamma)\gamma L(x) - \sigma(\gamma)\gamma  L p (x)\nonumber\\
   &=\sigma (\gamma) \gamma \left( (Dp)_{(x)}L - L p(x) \right)=\sigma (\gamma) \gamma \adl (p)(x),\nonumber
   \end{align}
so $\adl(p) \in \QQ_V(\Gamma)$. The other permutation is analogous, just use $\sigma^2(\gamma)=1,$ for all $\gamma \in \Gamma$.
\eproof

\begin{lemma} $\SS \adl = \adl \RR$.
\end{lemma}
\dem Note that
\begin{equation}
\label{eq:adlgammaestrela}
\adl(\gamma \star p)=\sigma(\gamma)\gamma \star (\adl(p)),~ \forall p \in \PP_V
 \end{equation}
is obtained by remarking that
\[  \adl(\gamma \star p)(x) = \adl(\gamma^{-1} p \gamma)(x),\] for all $x \in V.$ From this, we have
\begin{equation}
\label{eq:ss-adl}
\SS(\adl(p)) = \frac{1}{2} \left( \adl(p) + \adl(\delta \star p)   \right).
 \end{equation}
Now, just use the definitions of $\adl$ and $\RR$.
\eproof

We now present our main theorem. We have to find a complement to the homological operator $\adlk$ to the $\Gamma-$equivariants inside the vector space of $\Gamma-$reversible-equivariants. Our result is to recognize $\QQ_V^k(\s \rtimes \Gamma)$ as this complement. Consider the semidirect product $\s \rtimes \Gamma$ of the groups $\s$ and $\Gamma$ for which the homomorphism $\mu :\Gamma \to \aut(\s)$ is defined by
\begin{equation}\nonumber
\mu(\gamma)(e^{sL^t})=e^{\sigma(\gamma)sL^t}.
\end{equation}

\noindent By Proposition \ref{prop:representacao}, $\mu$ defines the action of $\s \rtimes \Gamma$ on $V$
\begin{equation}\nonumber
(e^{sL^t},\gamma)\cdot v = e^{sL^t}(\gamma v ).
\end{equation}
In this case, the equality (\ref{eq:acaoproduto}) is
\begin{equation}
\label{eq:acaonaocomuta}
e^{sL^t}(\gamma v )=\gamma  (e^{\sigma(\gamma)sL^t} v).
\end{equation}

\begin{theorem} \label{Teorema:principal} For $k \geq 2$, we have
\begin{equation}
\QQ_V^k(\Gamma)= \QQ_V^k(\s \rtimes \Gamma) \oplus \adlk(\PP_V^k(\Gamma)).
\label{eq1}
\end{equation}
\end{theorem}

To prove Theorem \ref{Teorema:principal}, we present three lemmas.  First, we use the action of $\Gamma$ on $\PP_V$, given in (\ref{eq:action-on-pol.-spaces}),  to define the mapping $\pi:\PP_V \to \PP_V$:
\begin{equation}
\label{eq:PROJECTION}
\pi(p)= \frac{1}{2} \left(\int_{\Gamma_+} \tau \star p ~d\tau - \int_{\Gamma_+}{(\delta \tau) \star p~d\tau}\right),
\end{equation}
where $\displaystyle\int_{\Gamma_+}$ is the normalized Haar integral over $\Gamma_{+}$. Notice that $\pi$ is an extension of the operator $\SS$. More then that, they have the same target space. In fact, we have:

\begin{lemma} \label{lemma1} The mapping $\pi:\PP_V \to \QQ_V(\Gamma)$ is a linear projection which preserves the grading of the algebra $\PP_V$.
\end{lemma}
\dem  By the linearity of the Haar integral it follows that $\pi$ is linear and  if $p$ has degree $k$, so does $\pi(p)$. To prove that $ \pi (p) \in \QQ_V^k(\Gamma)$, use
$$
\sigma(\gamma) \gamma \star \pi(p)= \frac{1}{2} \left( \int_{\Gamma_+} \sigma(\gamma)(\tau \gamma) \star p~d \tau - \int_{\Gamma_+}\sigma(\gamma) (\delta\tau\gamma) \star p~d\tau \right)
$$

\noindent and check that $\pi(p)=\sigma(\gamma)\gamma \star \pi(p)$ for elements $\gamma$ in $\Gamma_+$ and in $\delta \Gamma_+$ separately, using the left and right invariance of Haar integral and also the normality of the subgroup $\Gamma_+$. That $\pi^2=\pi $ follows from the fact that $\pi(p)=p$, $\forall p\in \QQ_V(\Gamma)$. \eproof

\begin{lemma} \label{lemma2} The projection $\pi$ satisfies  $\pi(\PP_V(\s))= \QQ_V(\s \rtimes \Gamma)$.
\end{lemma}
\dem
For $p\in \PP_V^k(\s)$ we want  $\pi(p) \in \QQ_V^k(\s \rtimes \Gamma)$,  that is, $\sigma(\gamma)(\gamma e^{\sigma(\gamma)sL^t})\star\pi(p) = \pi(p)$ for $\gamma \in \Gamma$. To prove this, we use (\ref{eq:acaonaocomuta}).  If $\gamma\in \Gamma_+$, we have
\begin{align}
\sigma(\gamma)(\gamma e^{\sigma(\gamma)sL^t}) \star  \pi(p) & = \frac{1}{2} \left( \int_{\Gamma_+}(\tau \gamma e^{sL^t}) \star p~d\tau - \int_{\Gamma_+}(\delta \tau \gamma e^{sL^t}) \star p~d \tau \right) \nonumber \\
& =\frac{1}{2} \left( \int_{\Gamma_+} (e^{sL^t}\tau ) \star p~d\tau - \int_{\Gamma_+}(e^{\sigma(\delta\tau)sL^t}\delta \tau ) \star p~d \tau \right) \nonumber \\
& =\frac{1}{2} \left( \int_{\Gamma_+} \tau \star ( e^{sL^t} \star p)~d\tau - \int_{\Gamma_+}(\delta \tau ) \star (e^{\sigma(\delta\tau)sL^t} \star p)~d \tau \right) \nonumber \\
& =\frac{1}{2} \left( \int_{\Gamma_+} \tau \star  p~d\tau - \int_{\Gamma_+} (\delta \tau) \star p~d\tau \right) =\pi(p), \nonumber
\end{align}
where the third equality follows from the definition of the action $\star$ given in (\ref{eq:action-on-pol.-spaces}) and the fourth equality follows because $p$ is $\s$-equivariant.

If $\gamma\in \delta\Gamma_+$  we have $\gamma=\delta \lambda$ for some $\lambda \in \Gamma_+$. Furthermore $\Gamma_+\triangleleft \Gamma$ so that $\delta \tau=\tilde{\tau}\delta $ for   $\tau, \tilde{\tau} \in \Gamma_+$. Then we obtain

\begin{align}
\sigma(\gamma)(\gamma e^{\sigma(\gamma)sL^t}) \star  \pi(p) & = \frac{1}{2}\left( \int_{\Gamma_+} (\delta\tau\gamma e^{-sL^t}) \star p~d\tau - \int_{\Gamma_+} (\tau \gamma e^{-sL^t}) \star p~ d\tau \right) \nonumber \\
&= \frac{1}{2}\left( \int_{\Gamma_+} (e^{-\sigma(\delta \tau \gamma)sL^t} \delta\tau\gamma)\star p ~d\tau - \int_{\Gamma_+}(e^{-\sigma(\tau\gamma)sL^t}\tau\gamma) \star p ~d\tau \right)\nonumber\\
&= \frac{1}{2}\left( \int_{\Gamma_+} ( \delta\tau\gamma)\star (e^{-sL^t} \star p) ~d\tau - \int_{\Gamma_+}(\tau\gamma)\star (e^{sL^t} \star p) ~d\tau \right)\nonumber\\
&= \frac{1}{2}\left( \int_{\Gamma_+} (\delta \tau \gamma)\star p ~d\tau - \int_{\Gamma_+} (\tau \gamma)\star p ~d \tau \right)\nonumber\\
& = \frac{1}{2}\left( \int_{\Gamma_+} ( \tilde{\tau} \delta^2\lambda) \star p~d \tilde{\tau} - \int_{\Gamma_+} (\delta \tilde{\tau} \lambda) \star p~d \tilde{\tau} \right) \nonumber \\
& = \frac{1}{2}\left( \int_{\Gamma_+}\tilde{\tau} \star p~ d\tilde{\tau }- \int_{\Gamma_+} (\delta \tilde{\tau}) \star p~d\tilde{\tau }\right)
=\pi(p). \nonumber
\end{align}

To prove the other inclusion, set $g\in \QQ_V^k(\s \rtimes \Gamma)$. By Proposition \ref{teoremainter}, $g \in \PP_V^k(\s)$ and $g=\sigma(\gamma)\gamma \star g$, for all $\gamma \in \Gamma$. Then,
$$
\pi(g)=\frac{1}{2} \left( \int_{\Gamma_+} \tau \star g ~d\tau - \int_{\Gamma_+}(\delta \tau) \star g~d\tau \right)
= \frac{1}{2} \left(\int_{\Gamma_+} g ~d\tau - \int_{\Gamma_+}  - g ~d\tau \right)
=  g. \eproof$$ 

For the next lemma, observe that
\begin{equation} \nonumber
\int_{\Gamma} \adl(\gamma\star p) ~d\gamma = \adl\int_{\Gamma}\gamma\star p~d\gamma,
\end{equation}
which follows directly from the linearity of $\adl$.

\begin{lemma} \label{lemma3} The projection $\pi$ satisfies $\pi(\adl (\PP_V))= \adl(\PP_V(\Gamma))$.
\end{lemma}

\dem  Let $p\in\PP_V$ and consider the equality (\ref{eq:adlgammaestrela}). Then,
\begin{align}
\pi(\adl (p)) & = \frac{1}{2} \left( \int_{\Gamma_+} \gamma \star (\adl (p)) d\gamma - \int_{\Gamma_+} (\delta \gamma)\star(\adl (p) )d\gamma \right) \nonumber \\
& =\frac{1}{2} \left(\int_{\Gamma_+} \sigma(\gamma) \adl(\gamma \star p)d\gamma - \int_{\Gamma_+} \sigma(\delta\gamma) \adl((\delta\gamma) \star p) d\gamma \right) \nonumber \\
& =\frac{1}{2} \left( \int_{\Gamma_+} \adl(\gamma \star p) d\gamma + \int_{\Gamma_+} \adl((\delta\gamma) \star p)d\gamma \right) \nonumber \\
& =\frac{1}{2} \left( \adl \int_{\Gamma_+} \gamma \star p~ d\gamma + \adl \int_{\Gamma_+} (\delta\gamma) \star p~d \gamma \right) \nonumber \\
& =\adl \left[ \frac{1}{2} \left( \int_{\Gamma_+} \gamma \star p~d\gamma + \int_{\Gamma_+} (\delta \gamma) \star p ~d\gamma \right) \right] \nonumber \\
& =\adl \left( \int_{\Gamma} \gamma \star p~ d\gamma \right). \nonumber
\end{align}

\noindent The last equality uses Fubini theorem (see \cite[Proposition I 5.16]{brocker}). Now, $\displaystyle\int_{\Gamma} \gamma \star p~ d\gamma \in \PP_V(\Gamma)$ and any element in $\PP_V(\Gamma)$ is of the form $\displaystyle\int_{\Gamma} \gamma \star p~ d\gamma$, for some $p \in \PP_V$.
\noindent \eproof

\demm
We apply the projection $\pi$ given in (\ref{eq:PROJECTION}) on equality (\ref{eq:ELPHICK}) and now we use Lemmas \ref{lemma1}, \ref{lemma2} and \ref{lemma3} to obtain
\begin{equation}
\QQ_V^k(\Gamma) = \QQ_V^k(\s \rtimes \Gamma) + \adlk (\PP_V^k(\Gamma)).
\label{eqdecomposicao3}
\end{equation}

By Proposition \ref{teoremainter}, $\QQ_V^k(\s \rtimes \Gamma)= \PP_V^k(\s) \cap \QQ_V^k(\Gamma)$, which together with (\ref{eq:ELPHICK}) gives
$$
\QQ_V^k(\s \rtimes \Gamma) \cap \adlk (\PP_V^k(\Gamma)) \subset \PP_V^k(\s) \cap \adl (\PP_V^k) = \{0\}.
$$ \eproof

From all the discussion of this section, the following result is now a direct consequence of Theorem \ref{Teorema:principal}.

\begin{theorem} \label{Teorema:principal2} Let $\Gamma$ be a compact Lie group acting linearly on $V$ and consider $X: V \to V$ a smooth $\Gamma-$reversible-equivariant vector field, $X(0)=0$ and $L=(dX)_0$. Then $(\ref{eq:sistemainicial})$ is formally conjugate to
\begin{equation} \nonumber
\dot{x}=Lx + g_2(x) + g_3(x) + \ldots
\end{equation}
where, for each $k \geq 2 $, $g_k \in \QQ_V^k(\s \rtimes \Gamma)$.
\end{theorem}

In practice we may be able to find the form of elements in $\QQ_V(\s\rtimes\Gamma)$. The main tool for that is given in \cite[Algorithm 3.7]{manoel1}. From that, one has just to select the general polynomial mapping in this module of the degree one wishes to truncate the normal form.


\section{Examples} \label{section:APPLICATIONS}

In this section we apply Theorems \ref{Teorema:principal} and \ref{Teorema:principal2} to deduce
the normal forms for two distinct examples: a  $\z_2-$reversible-equivariant vector field (with no
nontrivial symmetries) and a $\z_2 \times \z_2-$reversible-equivariant vector field.
Both of them are defined by a vector field on $\R^{6}$
\begin{equation}\label{eq:sistema-applications}
\dot{x}=X(x)
\end{equation}
whose linearization about the origin has matrix of type
\begin{equation}\label{eq:matrix-L}
L= \left(
\begin{matrix}            0 &1 & & & & \\
                          0 &0 & & & & \\
                                    &  &  0 &\omega_1  \\
                             & &  -\omega_1 & 0\\
                              &  &     &   & 0        & \omega_2 \\
                                &   &  &   &-\omega_2 & 0
\end{matrix}
\right),
\end{equation}
 with nonzero $\omega_1, \omega_2$ under a resonance condition $n_1 \omega_2 - n_2 \omega_1 = 0,$  $n_1, n_2 \in
\N$  nonzero. Under these conditions, the system (\ref{eq:sistema-applications}) is called $(n_1:n_2)-$resonant. The deduction of a normal form via $\ker \adltk $ becomes harder in computation as $n_1$ and $n_2$ get larger (we refer to \cite{lima} and \cite{mereu}, where the authors deal with some particular choices of $n_1$ and $n_2$). We emphasize that the usage of Theorem \ref{Teorema:principal2} circumvents  this problem, and works equality well for any values of $n_1$ and $n_2$. \\

By \cite[Proposition XVI 5.7]{golub},  in the present case we have $\s = \r \times \s^1,$
where $\r \cong\left\{\begin{pmatrix} 1 & 0 \\ s & 1 \end{pmatrix}, s\in \R \right\}$. Here we use complex coordinates and the diagonal action of $\s = \r \times \s^1$ on $\R^2 \times \C^2$, where
\begin{small}
\begin{equation}
\label{eq:action-R2xTn}
s(x_1,x_2)= (x_1, sx_1 + x_2) \quad \text{and} \quad \theta(z_1,z_2)= (e^{in_1\theta}z_1, e^{i n_2 \theta}z_2),
\end{equation}
\end{small} \noindent for $s \in \r$ and $\theta \in \s^1.$ For the action of $\r$ on $\R^2$, $\PP_{\small\R^2}(\r)$ is generated over  $\P_{\small\R^2}(\r)= \left< x_1 \right>$ by the set
\begin{equation}
\label{eq:GENERATORS}
\left\{ (x_1, x_2 ), ~( 0 , 1)  \right\}.
\end{equation}

Invariant and equivariant generators on $\C^2$ under $\s^1$ can be found in \cite[Theorem 4.2, Chapter XIX]{golub}:
$$
 |z_1|^2, |z_2|^2, \re(z_1^{n_2}\bar{z}_2^{n_1}), \im(z_1^{n_2}\bar{z}_2^{n_1})
$$
and
\begin{small}
$$
\begin{pmatrix} z_1 \\0  \end{pmatrix}, \begin{pmatrix} z_1 i \\0  \end{pmatrix}, \begin{pmatrix} \bar{z}_1^{n_2-1} z_2^{n_1} \\0  \end{pmatrix}, \begin{pmatrix} \bar{z}_1^{n_2-1} z_2^{n_1}i \\0  \end{pmatrix}, \begin{pmatrix} 0 \\ z_2 \end{pmatrix}, \begin{pmatrix} 0 \\ z_2i \end{pmatrix}, \begin{pmatrix} 0 \\ z_1^{n_2}\bar{z}_2^{n_1-1} \end{pmatrix} , \begin{pmatrix} 0 \\ z_1^{n_2}\bar{z}_2^{n_1-1}i \end{pmatrix},
$$
\end{small}
respectively.

\subsection{$\z_2-$reversible-equivariant normal form}

In this subsection we consider vector fields that are reversible equivariant under the action of $\z_2$ generated by
\begin{equation}
\label{eq:PHI}
\phi(x_1, x_2, z_1, z_2) = (x_1, -x_2, \bar{z_1}, \bar{z_2}),
\end{equation}
for $x_i \in \R$ and $z_i \in \C,$ with $i = 1,2.$ This involution is assumed to be a reversibility. The result is:

\begin{theorem} \label{teoremaformanormalressonante} Let $\dot{x}= Lx + h(x)$ a {\rm $\z_2-$}reversible-equivariant system, with $L$ defined in $(\ref{eq:matrix-L}).$ Then, this system is formally conjugate to:
\begin{small}
\begin{align}
\dot{x}_1 &= x_2 ~ +~  x_1\im(z_1^{n_2}\bar{z}_2^{n_1})f_{0}(X)\nonumber\\
\dot{x}_2 &= f_{1}(X)   ~+~    x_2\im(z_1^{n_2}\bar{z}_2^{n_1})f_{0}(X)\nonumber\\
\dot{z}_1 &= -~i\omega_1z_1 ~ +~  iz_1f_{2}(X) ~ +~   i\bar{z}_1^{n_2-1}z_2^{n_2}f_{3}(X) ~ + ~ z_1\im(z_1^{n_2}\bar{z}_2^{n_1})f_{4}(X)  \nonumber\\
          &~~~+ ~\bar{z}_1^{n_2-1}z_2^{n_2}\im(z_1^{n_2}\bar{z}_2^{n_1})f_{5}(X)\nonumber\\
\dot{z}_2 &= -~i\omega_2z_2 ~ +~  iz_2f_{6}(X) ~  +
~iz_1^{n_2}\bar{z}_2^{n_1-1}f_{7}(X)  ~+~
z_2\im(z_1^{n_2}\bar{z}_2^{n_1})f_{8}(X) \nonumber\\
          &~~~ +~ z_1^{n_2}\bar{z}_2^{n_1-1}\im(z_1^{n_2}\bar{z}_2^{n_1})f_{9}(X),
\nonumber
\end{align}
\end{small}

\noindent for some $f_i: \R^{4} \to \R$, $i= 0, \ldots, n$,
$X=(x_1,|z_1|^2, |z_2|^2, \re(z_1^{n_2}\bar{z}_2^{n_1}))$.
\end{theorem}

\dem Here the whole group is $\tilde{\Gamma} = (\r \times \s^1) \rtimes \z_2$. From Theorem \ref{Teorema:principal}, we need to compute the general form of elements in $\QQ_{\small{\R^2 \times \C^2}} ( \tilde\Gamma)$. We have that $\tilde{\Gamma}_+ = \r \times \s^1$.  From Lemma \ref{lemma:geradores},
$$
\{ x_1,  |z_1|^2, |z_2|^2, \re(z_1^{n_2}\bar{z}_2^{n_1}), \im(z_1^{n_2}\bar{z}_2^{n_1})\}
$$
is a Hilbert basis for $\P_{\small{\R^2 \times
\C^2}}(\tilde{\Gamma}_{+})$, and the generators for
$\PP_{\small{\R^2 \times \C^2}}(\tilde{\Gamma}_+)$ over the ring
$\P_{\small{\R^2 \times \C^2}}(\tilde{\Gamma}_+)$ are given by
$$
\begin{pmatrix} x_1 \\ x_2 \\0\\0 \end{pmatrix}, \begin{pmatrix} 0 \\ 1 \\0\\0\end{pmatrix},
\begin{pmatrix} 0\\0\\ z_1\\  0  \end{pmatrix}, \begin{pmatrix} 0\\0\\iz_1\\  0  \end{pmatrix}, \begin{pmatrix} 0\\0\\\bar{z}_1^{n_2-1} z_2^{n_1} \\0  \end{pmatrix}, \begin{pmatrix}, 0\\0\\i\bar{z}_1^{n_2-1} z_2^{n_1} \\0  \end{pmatrix},$$ $$ \begin{pmatrix} 0\\0\\0 \\ z_2 \end{pmatrix}, \begin{pmatrix} 0\\0\\0 \\ iz_2 \end{pmatrix}, \begin{pmatrix} 0\\0\\0 \\ z_1^{n_2}\bar{z}_2^{n_1-1} \end{pmatrix}, \begin{pmatrix} 0\\0\\0 \\ iz_1^{n_2}\bar{z}_2^{n_1-1} \end{pmatrix}.
$$

Then, we apply the algorithm \cite[Algorithm 3.7]{manoel1} to obtain generators for
$\QQ_{\small{\R^2 \times \C^2}}(\tilde{\Gamma})$ over $\P_{\small{\R^2 \times \C^2}}(\tilde{\Gamma})$:
\begin{small}
$$
\begin{pmatrix} x_1\im(z_1^{n_2}\bar{z}_2^{n_1}) \\ x_2\im(z_1^{n_2}\bar{z}_2^{n_1}) \\0\\0 \end{pmatrix},
\begin{pmatrix} 0\\0\\iz_1\\  0  \end{pmatrix},  \begin{pmatrix} 0\\0\\ i\bar{z}_1^{n_2-1} z_2^{n_1} \\0  \end{pmatrix},
\begin{pmatrix} 0\\0\\ z_1 \im(z_1^{n_2}\bar{z}_2^{n_1}) \\0  \end{pmatrix},
\begin{pmatrix} 0\\0\\ \bar{z}_1^{n_2-1} z_2^{n_1}\im(z_1^{n_2}\bar{z}_2^{n_1}) \\0  \end{pmatrix},
$$
$$
\begin{pmatrix} 0\\1\\0\\0 \end{pmatrix},\begin{pmatrix} 0\\0\\0 \\ iz_2 \end{pmatrix}, \begin{pmatrix} 0\\0\\0 \\ iz_1^{n_2}\bar{z}_2^{n_1-1} \end{pmatrix},
\begin{pmatrix} 0\\0\\0 \\ z_2\im(z_1^{n_2}\bar{z}_2^{n_1}) \end{pmatrix},
\begin{pmatrix} 0\\0\\0 \\ z_1^{n_2}\bar{z}_2^{n_1-1} \im(z_1^{n_2}\bar{z}_2^{n_1}) \end{pmatrix}.
$$
\end{small}

Using the Reynolds operators, $R$ and $S$, and Theorem \ref{teoremagamma}, we obtain
\begin{equation} \label{eq:base-hilbert-ressonante}
\{x_1, |z_1|^2, |z_2|^2, \re(z_1^{n_2}\bar{z}_2^{n_1})\}
\end{equation}
as a Hilbert basis for $\P_{\small{\R^2 \times \C^2}} (\tilde{\Gamma})$.

In fact, we have $R(x_1) = x_1,$ $R(|z_1|^2) = |z_1|^2,$ $R(|z_2|^2)
= |z_2|^2,$ $R(\re(z_1^{n_2}\bar{z}_2^{n_1})) =
\re(z_1^{n_2}\bar{z}_2^{n_1})$, $R(\im(z_1^{n_2}\bar{z}_2^{n_1})) =
0$, $S(x_1)$ $=$ $S(|z_1|^2) =$ $S(|z_2|^2) =$
$S(\re(z_1^{n_2}\bar{z}_2^{n_1})) = 0$ and
$S(\im(z_1^{n_2}\bar{z}_2^{n_1})) = \im(z_1^{n_2}\bar{z}_2^{n_1})$.
But $(\im(z_1^{n_2}\bar{z}_2^{n_1}))^2$ is obtained from
(\ref{eq:base-hilbert-ressonante}).\eproof

\subsection{$\z_2 \times \z_2-$reversible-equivariant normal form}

In this subsection we consider vector fields that are reversible equivariant under the action of $\z_2^\phi \times \z_2^\psi$ generated by the involutions $\phi$ as in (\ref{eq:PHI}) and

\begin{equation}
\label{eq:PSI}
\psi(x_1, x_2, z_1, z_2)= (a_0x_1, -a_0x_2, a_1\bar{z_1}, a_2\bar{z_2}),
\end{equation}
with $a_i = \pm 1,$ for $i = 0,1,2.$ These involutions are assumed to be reversibilities.

We want to obtain a $\z_2^\phi \times \z_2^\psi-$reversible-equivariant normal form of (\ref{eq:sistema-applications}). In this case, the whole group is $\tilde{\Gamma} = (\r \times \s^1) \rtimes (\z_2^\phi \rtimes \z_2^\psi)$ and from Theorem \ref{Teorema:principal}, we need to compute the generators for $\QQ_{\small{\R^2 \times \C^2}} ( \tilde\Gamma)$ under $\P_{\small{\R^2 \times \C^2}} (\tilde{\Gamma})$. We have that $\tilde{\Gamma}_+ = (\r \times \s^1) \rtimes \z_2^{\phi \psi}$ and apply
the algorithm \cite[Algorithm 3.7]{manoel1} to obtain the result.

The possible values of $a_0, a_1$ e $a_2$ in (\ref{eq:PSI}) give rise to four types of normal forms which are described in Table \ref{tab:tabletypes}, with respective generators given in Table \ref{tab:tablegenerators}. Here, $u_1(x,z)= x_1,$ $u_2(x,z)= |z_1|^2,$ $u_3(x,z)=|z_2|^2,$ $u_4(x,z)= \re(z_1^{n_2}\bar{z}_2^{n_1})$ and $u_5(x,z)=\im(z_1^{n_2}\bar{z}_2^{n_1}).$ Also,

\begin{itemize}
\item[] $H_0(x,z)=\begin{pmatrix} 0,1,0,0 \end{pmatrix}$,

\item[] $H_1(x,z)=\begin{pmatrix} x_1\im(z_1^{n_2}\bar{z}_2^{n_1}) , x_2\im(z_1^{n_2}\bar{z}_2^{n_1}) ,0,0 \end{pmatrix}$,

\item[] $H_2(x,z)= \begin{pmatrix} 0,0,iz_1,  0  \end{pmatrix}$,

\item[] $H_3(x,z)= \begin{pmatrix} 0,0, i\bar{z}_1^{n_2-1} z_2^{n_1} ,0  \end{pmatrix}$,

\item[] $H_4(x,z)=\begin{pmatrix} 0,0, z_1 \im(z_1^{n_2}\bar{z}_2^{n_1}) ,0  \end{pmatrix}$,

\item[] $H_5(x,z)=\begin{pmatrix} 0,0, \bar{z}_1^{n_2-1} z_2^{n_1}\im(z_1^{n_2}\bar{z}_2^{n_1}) ,0  \end{pmatrix}$,

\item[] $H_6(x,z)=\begin{pmatrix} 0,0,0 , iz_2 \end{pmatrix}$,

\item[] $H_7(x,z)=\begin{pmatrix} 0,0,0, iz_1^{n_2}\bar{z}_2^{n_1-1} \end{pmatrix}$,

\item[] $H_8(x,z)=\begin{pmatrix} 0,0,0, z_2\im(z_1^{n_2}\bar{z}_2^{n_1}) \end{pmatrix}$,

\item[] $H_9(x,z)=\begin{pmatrix} 0,0,0, z_1^{n_2}\bar{z}_2^{n_1-1} \im(z_1^{n_2}\bar{z}_2^{n_1}) \end{pmatrix}$.
\end{itemize}

\begin{table}
    \centering
\begin{tabular}{|c|c|c|c|}
\hline
$a_0$ & $a_1, \  a_2$ & $n_1, \ n_2$ & Type \\
\hline
  $a_0=1$ & $a_1=a_2=1 $ &  --- & A \\
          \cline{2-4}
          &  $a_1=a_2=-1$ & $n_1 + n_2$ even &   A\\
          \cline{3-4}
          &             & $n_1 + n_2$ odd &   B\\
          \cline{2-4}
          &  $ a_1=-a_2=1$  &  $n_1$ even  &   A\\
          \cline{3-4}
          &                 &  $n_1$ odd &   B\\
          \cline{2-4}
          &  $a_1=-a_2 =-1$ &  $n_2$ even &     A\\
          \cline{3-4}
          &                 &  $n_2$ odd &   B\\
                       \hline
        $a_0=-1$ & $a_1=a_2=1 $ &  --- &  C \\
        \cline{2-4}
        &  $a_1=a_2=-1$ & $n_1 + n_2$ even &   C\\
        \cline{3-4}
        &             & $n_1 + n_2$ odd &   D\\
        \cline{2-4}
        &  $ a_1=-a_2=1$  &  $n_1$ even  &   C\\
        \cline{3-4}
       &                 &  $n_1$ odd &  D\\
       \cline{2-4}
       &  $a_1=-a_2=-1$ &  $n_2$ even &     C\\
       \cline{3-4}
       &                 &  $n_2$ odd &   D\\
       \hline
      \end{tabular}
\caption{Normal forms of $\z_2 \times \z_2-$reversible-equivariant vector fields}
    \label{tab:tabletypes}
\end{table}

\begin{table}
    \centering
\begin{tabular}{|c|c|c|}
\hline
Type & $\QQ_{\small\R^2 \times \C^2}(\tilde{\Gamma})$ & $\P_{\small\R^2 \times \C^2}(\tilde{\Gamma})$\\
\hline
 A &   $H_j$, $~0 \leq j \leq 9$     &   $u_1, u_2, u_3, u_4$\\
\hline
 B  &   $H_k$, $~u_4 H_l~$ for $k=0,2,5,6,9 $ and $l=1,3,4,7,8$    &  $u_1,u_2,u_3,u_4^2$ \\
\hline
 C  &  $u_1H_0$, $~H_k~$ for $1 \leq k \leq 9$     & $u_1^2, u_2,u_3,u_4$\\
\hline
 D  &  $H_k$, $~u_1H_l$,  $~u_4 H_l~$         &  $u_1^2, u_2,u_3,u_4^2,u_1u_4$    \\
        &      for $k=2,5,6,9$ and $l=0,1,3,4,7,8$ &  $u_1^2, u_2,u_3,u_4^2,u_1u_4$    \\
\hline
\end{tabular}
    \caption{Generators in $\R^2 \times \C^2$ for each possible type given in Table~{tab:tabletypes}}
    \label{tab:tablegenerators}
\end{table}

To determine the Hilbert basis for $\P_{\small{\R^2 \times \C^2}}(\tilde{\Gamma})$ we use the Theorem \ref{teoremagamma}, by considering the Reynolds operators $R$, $S$ and the Hilbert basis for $\P_{\small{\R^2 \times \C^2}}((\r \times \s^1) \rtimes \z_2^\phi).$ The result is:

\begin{theorem} Let $\dot{x}= X(x)$ a {\rm $\z_2^\phi \times \z_2^\psi-$}reversible-equivariant system, with $L=(dX)_0$ defined in $\eqref{eq:matrix-L}.$ Then, this system is formally conjugate to one of the following systems, types given according to Table~\ref{tab:tabletypes}:

\noindent Type A:
\begin{small}
\begin{align}
\dot{x}_1 &= x_2  ~+~  x_1\im(z_1^{n_2}\bar{z}_2^{n_1})f_{0}(X)\nonumber\\
\dot{x}_2 &= f_{1}(X)  ~ +  ~  x_2\im(z_1^{n_2}\bar{z}_2^{n_1})f_{0}(X)\nonumber\\
\dot{z}_1 &= -~i\omega_1z_1  ~+~  iz_1f_{2}(X) ~ + ~  i\bar{z}_1^{n_2-1}z_2^{n_2}f_{3}(X) ~ +~  z_1\im(z_1^{n_2}\bar{z}_2^{n_1})f_{4}(X)\nonumber\\
          &~~~  +~ \bar{z}_1^{n_2-1}z_2^{n_2}\im(z_1^{n_2}\bar{z}_2^{n_1})f_{5}(X)\nonumber\\
\dot{z}_2 &= -~i\omega_2z_2 ~ +~  iz_2f_{6}(X)  ~ +~iz_1^{n_2}\bar{z}_2^{n_1-1}f_{7}(X)~+~z_2\im(z_1^{n_2}\bar{z}_2^{n_1})f_{8}(X) \nonumber\\
          &~~~+ z_1^{n_2}\bar{z}_2^{n_1-1}\im(z_1^{n_2}\bar{z}_2^{n_1})f_{9}(X)
\nonumber
\end{align}
\end{small}
\noindent for  $f_i: \R^{4},0 \to \R$, $0 \leq i \leq 9$,  $X=(x_1, |z_1|^2, |z_2|^2, \re(z_1^{n_2}\bar{z}_2^{n_1}))$.\\

\noindent Type B:
\begin{small}
\begin{align}
\dot{x}_1 &= x_2 ~ + ~ x_1 \re(z_1^{n_2}\bar{z}_2^{n_1}) \im(z_1^{n_2}\bar{z}_2^{n_1})f_{0}(X)\nonumber\\
\dot{x}_2 &= f_{1}(X)  ~ +  ~  x_2 \re(z_1^{n_2}\bar{z}_2^{n_1}) \im(z_1^{n_2}\bar{z}_2^{n_1})f_{0}(X)\nonumber\\
\dot{z}_1 &= -~i\omega_1z_1 ~ +  ~iz_1f_{2}(X) ~ +~  \bar{z}_1^{n_2-1}z_2^{n_2}\im(z_1^{n_2}\bar{z}_2^{n_1})f_{3}(X) \nonumber\\
          &~~~ + ~ i\bar{z}_1^{n_2-1}z_2^{n_2}\re(z_1^{n_2}\bar{z}_2^{n_1})f_{4}(X) ~ + ~ z_1\re(z_1^{n_2}\bar{z}_2^{n_1})\im(z_1^{n_2}\bar{z}_2^{n_1})f_{5}(X) \nonumber\\
\dot{z}_2 &= -~i\omega_2z_2 ~ +~  iz_2f_{6}(X) ~ + ~z_1^{n_2}\bar{z}_2^{n_1-1}\im(z_1^{n_2}\bar{z}_2^{n_1})f_{7}(X) \nonumber\\
          &~~~+ ~ i z_1^{n_2}\bar{z}_2^{n_1-1} \re(z_1^{n_2}\bar{z}_2^{n_1})f_{8}(X) ~ + ~z_2\re(z_1^{n_2}\bar{z}_2^{n_1})\im(z_1^{n_2}\bar{z}_2^{n_1})f_{9}(X)   \nonumber
\end{align}
\end{small}
\noindent for  $f_i: \R^{4},0 \to \R$, $0 \leq i \leq 9$,   $X=(x_1, |z_1|^2, |z_2|^2, \re^2(z_1^{n_2}\bar{z}_2^{n_1}))$.\\

\noindent Type C:
\begin{small}
\begin{align}
\dot{x}_1 &= x_2  ~+ ~ x_1\im(z_1^{n_2}\bar{z}_2^{n_1})f_{0}(X)\nonumber\\
\dot{x}_2 &= x_1f_{1}(X) ~  + ~   x_2\im(z_1^{n_2}\bar{z}_2^{n_1})f_{0}(X)\nonumber\\
\dot{z}_1 &= -~i\omega_1z_1  ~+~  iz_1f_{2}(X) ~ +~   i\bar{z}_1^{n_2-1}z_2^{n_2}f_{3}(X) ~ +~  z_1\im(z_1^{n_2}\bar{z}_2^{n_1})f_{4}(X)  \nonumber\\
          &~~~+ \bar{z}_1^{n_2-1}z_2^{n_2}\im(z_1^{n_2}\bar{z}_2^{n_1})f_{5}(X)\nonumber\\
\dot{z}_2&=-~i\omega_2z_2~+~iz_2f_{6}(X)~+~iz_1^{n_2}\bar{z}_2^{n_1-1}f_{7}(X) ~+ ~z_2\im(z_1^{n_2}\bar{z}_2^{n_1})f_{8}(X)\nonumber\\
         &~~~+z_1^{n_2}\bar{z}_2^{n_1-1}\im(z_1^{n_2}\bar{z}_2^{n_1})f_{9}(X)
\nonumber
\end{align}
\end{small}
\noindent for  $f_i: \R^{4},0 \to \R$, $0 \leq i \leq 9$,   $X=(x_1^2, |z_1|^2, |z_2|^2, \re(z_1^{n_2}\bar{z}_2^{n_1}))$.\\

\noindent Type D:
\begin{small}
\begin{align}
\dot{x}_1 &= x_2  ~+~  x_1^2\im(z_1^{n_2}\bar{z}_2^{n_1})f_{0}(X) ~ +~   x_1 \re(z_1^{n_2}\bar{z}_2^{n_1})\im(z_1^{n_2}\bar{z}_2^{n_1}) f_{1}(X) \nonumber\\
\dot{x}_2 &= x_1f_{2}(X)  ~ + ~   x_1 x_2\im(z_1^{n_2}\bar{z}_2^{n_1})f_{0}(X) ~+~  \re(z_1^{n_2}\bar{z}_2^{n_1}))f_{3}(X)    \nonumber\\
          &~~~+  x_2\re(z_1^{n_2}\bar{z}_2^{n_1})\im(z_1^{n_2}\bar{z}_2^{n_1})f_{1}(X)    \nonumber\\
\dot{z}_1 &= -~i\omega_1z_1  ~+~  iz_1f_{4}(X)~  +~\bar{z}_1^{n_2-1}z_2^{n_2}\im(z_1^{n_2}\bar{z}_2^{n_1})f_{5}(X) ~+~   ix_1\bar{z}_1^{n_2-1}z_2^{n_2}f_{6}(X)  \nonumber \\
          & ~~~+~ x_1 z_1\im(z_1^{n_2}\bar{z}_2^{n_1})f_{7}(X)~ + ~i\bar{z}_1^{n_2-1}z_2^{n_2} \re(z_1^{n_2}\bar{z}_2^{n_1})f_{8}(X)  \nonumber \\
          & ~~~+~z_1 \re(z_1^{n_2}\bar{z}_2^{n_1}) \im(z_1^{n_2}\bar{z}_2^{n_1})f_{9}(X)  \nonumber\\
\dot{z}_2&=-~i\omega_2z_2 ~ + ~iz_2f_{10}(X)+z_1^{n_2}\bar{z}_2^{n_1-1}\im(z_1^{n_2}\bar{z}_2^{n_1})f_{11}(X) ~+~ix_1z_1^{n_2}\bar{z}_2^{n_1-1}f_{12}(X)\nonumber\\
         & ~~~+ ~x_1 z_2\im(z_1^{n_2}\bar{z}_2^{n_1})f_{13}(X) ~+~ iz_1^{n_2}\bar{z}_2^{n_1-1}\re(z_1^{n_2}\bar{z}_2^{n_1})f_{14}(X)\nonumber\\
         &~~~+~ z_2 \re(z_1^{n_2}\bar{z}_2^{n_1})\im(z_1^{n_2}\bar{z}_2^{n_1})f_{15}(X) \nonumber
\end{align}
\end{small} \noindent for  $f_i: \R^{4},0 \to \R$, $0 \leq i \leq
15$, $X=(x_1^2, |z_1|^2, |z_2|^2, \re^2(z_1^{n_2}\bar{z}_2^{n_1}),
x_1\re(z_1^{n_2}\bar{z}_2^{n_1}))$.
\end{theorem}

\begin{remark}
The case $\phi = \psi$ is the $\z_2$-reversible case. In \cite{lima}
the authors present a truncated normal form for this case,
truncation at a low degree. Without resonance, the authors also
present the formal normal form, computed by solving the PDEs that
define the kernel of ${\adltk},$ for $k \geq 2$.
\end{remark}

\noindent {\it Acknowledgments.} The research of M.M. was supported by FAPESP,  BPE grant 2013/11108-7.


\end{document}